\documentclass[12pt]{article}
\usepackage{amssymb,amsmath,eucal,amsthm}
\numberwithin{equation}{section}

\newcommand{\R}{\mathbb{R}}
\newcommand{\Z}{\mathbb{Z}}
\newcommand{\C}{\mathbb{C}}
\newcommand{\PP}{\mathbb{P}}
\newcommand{\Q}{\mathbb{Q}}
\newcommand{\T}{\mathbb{T}}
\newcommand{\bA}{\mathbb{A}}
\newcommand{\K}{\mathbb{K}}
\newcommand{\OO}{\mathbb{O}}

\newcommand{\fh}{\mathfrak{h}}
\newcommand{\fv}{\mathfrak{v}}

\newcommand{\la}{\lambda}
\newcommand{\al}{\alpha}

\newcommand{\q}{/\!\!/}
\newcommand{\hf}{h}

\DeclareMathOperator{\spec}{spec}
\DeclareMathOperator{\Hom}{Hom}
\DeclareMathOperator{\Lie}{Lie}
\DeclareMathOperator{\Vol}{Vol}
\DeclareMathOperator{\vol}{vol}
\DeclareMathOperator{\Pic}{Pic}
\DeclareMathOperator{\tr}{tr}

\DeclareMathOperator{\ord}{ord}
\DeclareMathOperator{\codim}{codim}
\DeclareMathOperator{\supp}{supp}

\newtheorem{theorem}{Theorem}
\newtheorem{conjecture}{Conjecture}

\theoremstyle{definition}
\newtheorem{definition}{Definition}

\begin{document}

\title{Why would multiplicities be log-concave ?}
\author{Andrei Okounkov}
\date{}

\maketitle 

\begin{abstract}
It is a basic property of the entropy in statistical physics that is
concave as a function of energy. The analog of this in
representation theory would be the concavity of the logarithm of the multiplicity
of an irreducible representation as a function of its highest weight.
We discuss various situations where such concavity can be established or
reasonably conjectured and consider some implications of this concavity.
These are rather informal notes based on a number of talks I gave on the
subject, in particular, at the 1997 International Press lectures at UC Irvine. 
\end{abstract}

\setcounter{section}{-1}
\section{Introduction} 

The aim of these notes is to discuss some heuristic arguments, conjectures,
and rigorous results related to the following phenomenon. Physical analogy,
explained in Section \ref{s1}, suggests that under certain circumstances 
the logarithms of multiplicities of irreducible representations can be
expected to be concave as a function of the highest weight. In Section \ref{s2}
we discuss some cases when this is known or expected to be the case and
explore various implications of this concavity. In Section \ref{s3} we discuss
the classical limit, in which much more general results can be established.

This text is not a survey. It is based on several talks I gave on the subject 
on various occasions and represents only my personal point of view. I hope that
the somewhat informal style of these notes will make the basic ideas easier
to explain. For missing details, the reader is referred to the original papers
\cite{Ok1,Gr,Ok2,Ka,Ok3}. For surveys on log-concavity in general, see \cite{Br,St}. 

I very much benefited from the discussions with a number of people, first
of all, with my colleagues
from the Institute of Problems of Information Transmission, especially
R.~Dobrushin, G.~Olshanski, and S.~Pirogov, and also with V.~Ginzburg, W.~Graham, 
A.~Khovanskii, and A.~Kirillov. In particular, the results of \cite{Ok1} lead V.~Ginzburg
to conjecture that the push-forward of the Liouville measure on an arbitrary
symplectic manifold under the moment map for a compact group action 
should be log-concave. Same conjecture, independently of \cite{Ok1}, was proposed
by A.~Knutson (later, a counterexample to this conjecture was found 
in \cite{Ka}; for positive results see \cite{Gr,Ok2,Ok3}).  

I would like to thank 
A.~Buch for providing me with a program for computation of tensor product multiplicities.

\section{Physical motivation: entropy and its concavity}\label{s1} 

\subsection{}
Consider a quantum mechanical system, that is, 
a selfadjoint operator $H$ in a Hilbert space $V$. For simplicity, we
assume that $V$ is spanned by the eigenvectors of $H$. 

The multiplicity 
$\Omega(E)$ of an eigenvalue $E\in\spec(H)$ measures how many states of our system have 
the energy $E$. In other words, fixing an energy level $E$ this does not determine 
the state of the system uniquely: there remain
$\Omega(E)$ possibilities. The size of this indeterminacy equals $\log_2 \Omega(E)$
bits of information.  

In statistical physics, there is the basic relation\footnote
{This relation, in the form $S=k\log W$, is written on Boltzmann's tombstone.}
$$
\fbox{$\displaystyle S=k \log \Omega$}
$$
where $\Omega$ is the number
of states with given values of macroscopic parameters
such as energy, $k$ is the Boltzmann constant, and $S$ is the \emph{entropy},
which measures the degree of disorder in the system or, in other words, the
lack of information about the precise state of the system. We are thus led
to think of 
$$
S(E)=\log\Omega(E)
$$ 
as of the entropy of the energy level $E$.

\subsection{}

In statistical mechanics, the entropy is always a concave function of all additive 
macroscopic parameters such as energy $E$, volume $V$, or the
number of particles $N$.  
There
is a simple physical argument for this concavity and it goes as follows. 
Suppose we have two systems with parameters $(E_1,V_1,\dots)$ and
$(E_2,V_2,\dots)$, respectively, contained in two reservoirs 
separated by an impervious wall:
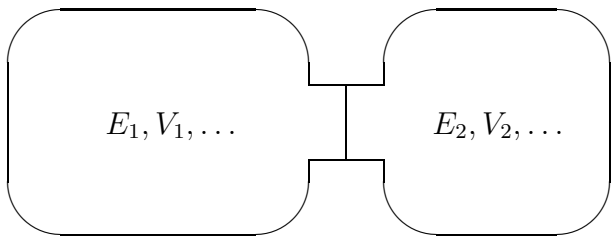
\begin{center}
\setlength{\unitlength}{1 cm}
\begin{picture}(8,5)
\put(2,3){\oval(4,2)[t]}
\put(2,2){\oval(4,2)[b]}
\put(1.3,2.35){$E_1,V_1,\dots$}
\put(0,2){\line(0,1){1}}
\put(4,2){\line(1,0){1}}
\put(4,3){\line(1,0){1}}
\put(4.5,2){\line(0,1){1}}
\put(6.5,2){\oval(3,2)[b]}
\put(6.5,3){\oval(3,2)[t]}
\put(8,2){\line(0,1){1}}
\put(5.65,2.35){$E_2,V_2,\dots$}
\end{picture}
\end{center}
Let us now bring them in contact by removing this wall. The  energy and the
volume of the new system will be $E_1+E_2$ and $V_1+V_2$, respectively, whereas
the entropy will increase
\begin{equation}\label{entr_mix}
S(E_1+E_2,\dots) \ge S(E_1,\dots) + S(E_2,\dots)
\end{equation}
because of the additional disorder introduced by allowing the systems to mix. 
The net increase in entropy is called the \emph{entropy of mixing} and  its 
positivity reflects the irreversibility of mixing. 

There is, however, one case when the mixing is clearly reversible and that is
when the two systems were identical to begin with, that is, when
$$
(E_1,V_1,\dots)=(E_2,V_2,\dots) \,,
$$
in which case we can simply insert back the wall to recover the original
situation. Thus, in this case the entropy of mixing vanishes. 
In other words, 
\begin{equation}\label{gibbs_p}
S(2E,\dots)=2S(E,\dots)\,.
\end{equation}
Combining \eqref{entr_mix} with \eqref{gibbs_p} we get the concavity of the entropy\footnote{
In thermodynamics, one has the relation $\frac{\partial S}{\partial E}=\frac1T$ where $T$ is
the temperature. Therefore $\frac{\partial^2 S}{\partial E^2}<0$ means that temperature rises when energy
increases.}.  

\subsection{}\label{s_irr}

Of course, in order to apply statistical considerations one needs
the system in question to have a very large or
infinite number of degrees of freedom.  
Still, it is natural to ask whether in some interesting cases one can
expect or, even better, prove the concavity of $S(E)$.  Fortunately, 
interesting examples do exists.

An obvious limitation for the entropy concavity principle is that the concavity
of $S(E)$ is clearly not preserved under direct sums. Hence our system has
be in some sense \emph{irreducible}. The concrete meaning of this irreducibility
will be different in different context. In Section \ref{s2}, the space $V$ will
be an irreducible module of some ambient group. In Section \ref{s3}, we will be dealing
with group actions on irreducible algebraic varieties.

\subsection{}

First, however, one 
has to modify the definition of concavity. Indeed, the support $\spec(H)$
of the function $\Omega(E)$ is countable and hence $S(E)$ cannot be concave in the
usual sense. 

We suppose that $\spec(H)$ is contained in a lattice, which without loss
of generality we can take to be $\Z\subset\R$, and we define concavity to mean 
$$
S(\al E_1+ (1-\al) E_2) \ge \al S(E_1)+ (1-\al) S(E_2) \,, \quad 
\al\in[0,1]\,,
$$
whenever the middle point $\al E_1+ (1-\al) E_2$ lies in the lattice $\Z$.  

The abstract form of this convention is the following: 

\begin{definition} Let $F:\bA\to\OO$ be a function from a Abelian
semigroup $\bA$ to an ordered Abelian 
semigroup $\OO$. We say that this function is \emph{concave} if
$$
(p+q) \,F(C) \ge p\, F(A) + q \,F(B)
$$
for any $A,B,C\in\bA$ satisfying
$$
(p+q)\, C = p\,  A+  q \, B \,, \quad p,q\in\Z_{\ge 0} \,.
$$
\end{definition}

In the case when
$$
\OO=(\R_{\ge 0},\times)
$$
is the multiplicative semigroup
of nonnegative real numbers with the usual ordering, we 
also call the function $F$ logarithmically concave, or
\emph{log-concave} for short. 

In our examples, $\bA$ will be usually isomorphic to $\Z^n$ or $\R^n$,
whereas the target semigroup $\OO$ will occasionally be something
more interesting. 

\subsection{}\label{s1K}

Since the eigenvalues of $H$ are now integers, the time evolution $e^{itH}$ 
defines a representation of the standard circle $\T^1$ on $V$.

More generally, for any compact group\footnote{
Here and in what follows we assume all compact groups to be \emph{connected}.}
 $K$, one can ask 
whether for some interesting representation $V$ of $K$ 
the multiplicities $\Omega(\la)$ of irreducible representations $V^\la$ 
$$
V = \bigoplus_{\la\in K^\wedge} \Omega(\la)\, V^\la
$$
form a log-concave function on the weight lattice of $K^\wedge$.

Examples of such representations will be discussed in Section \ref{s2}.
They are, in a sense, related to the 
``thermodynamics'' of classical groups.

\subsection{}\label{s11}

Now consider a Hamiltonian system of classical mechanics, that is, a
manifold $M^{2n}$ with symplectic form $\omega$ and with an energy function
$$
\hf: M^{2n}\to \R\,.
$$
The form $\dfrac{\omega^n}{n!}$ is a volume form on $M^{2n}$ which defines 
a measure (called the \emph{Liouville measure}). Let $\Omega(E)$ be the density
of  the push-forward of this measure under $\hf$
$$
\hf_*\left(\dfrac{\omega^n}{n!}\right) = \Omega(E)\, dE \,, \quad E\in\R \,.
$$
In other words, $\Omega(E)$ tells us 
how many states of our system have the energy $E$. Again, we think of 
$$
S(E)=\log S(E)
$$
as of the entropy\footnote{
This entropy is not to be confused with the entropy of the dynamical system defined
on $\hf^{-1}(E)$ by the Hamiltonian flow $\dot x=\{\hf,x\}$.}
of the energy level $E$. 

As in Section \ref{s1K}, this can be generalized to the situation of a 
Hamiltonian action of a compact group $K$ action on $M^{2n}$. Let
$$
\phi: M^{2n} \to \Lie(K)^*
$$ 
be the moment map for this action. For any $\xi\in  \Lie(K)^*$ the volume 
$$
\Omega(\xi)=\Vol\phi^{-1}(\xi)
$$
measures how many points of $X$ have the energy
$\xi$. This function is clearly
invariant under the coadjoint action of $K$ on $\Lie(K)^*$, so we can
and will assume that $\xi$ lies in the positive Weyl chamber $\fh_+$.

We can ask whether for some actions the function
$\log \Omega(\xi)$ is concave on $\fh_+$. Observe that such a concavity
implies, in particular, that the set 
$$
\supp\Omega(\xi)=\phi(M^{2n})\cap\fh_+\,,
$$
is convex, which is a famous classical result \cite{A,GS2,GS3,Ki}. 

It turns out that the supply of cases where $\log \Omega(\xi)$ is concave 
is now much richer than in the quantum situation. As shown by W.~Graham in \cite{Gr}, it includes
all torus actions on compact K\"ahler manifolds. It also includes \cite{Ok2}
all actions on projective varieties, possibly singular. It was
conjectured by V.~Ginzburg and A.~Knutson that it is true for any
symplectic $M^{2n}$. This was shown to be not the case by Y.~Karshon 
in \cite{Ka}. 

We will discuss this classical situation in algebraic setting in Section \ref{s3}. 

\section{Some results and conjectures on
log-concavity of multiplicities}\label{s2} 

\subsection{} 
Again, we begin with a  motivation, this time a historical one. 
Here is how the question of logarithmic concavity of multiplicities arose in the
``thermodynamics'' of classical groups. 

Let $U(\infty)$ denote the inductive limit of $U(n)$ with respect to 
standard embeddings $U(n)\subset U(n+1)$ which can be visualized as follows:
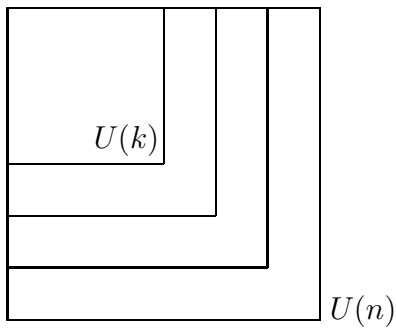
\begin{center}
\setlength{\unitlength}{0.69 cm}
\begin{picture}(6,6)
\put(0,0){\framebox(6,6)}
\put(0,1){\line(1,0){5}}
\put(5,1){\line(0,1){5}}
\put(0,2){\line(1,0){4}}
\put(4,2){\line(0,1){4}}
\put(0,3){\line(1,0){3}}
\put(3,3){\line(0,1){3}}
\put(6.2,0){$U(n)$}
\put(1.65,3.3){$U(k)$}
\end{picture}
\end{center}
The description of the characters\footnote{
An abstract definition of characters
is: indecomposable central continuous positive
definite functions.
 More concretely, they are spherical functions of the Gelfand pair
$$
U(\infty)\times U(\infty)\supset \textup{diag}\,U(\infty)
$$
or, equivalently, traces of
factor representations of type $\textup{I}_n$ or $\textup{II}_1$.} 
of $U(\infty)$ is a fundamental result with nontrivial history. 
Voiculescu in \cite{Vo} proved that 
functions of the form
\begin{equation}\label{char}
g\mapsto \det \left[e^{\gamma^+(g-1)+\gamma^-(g-1)}\, 
\prod \frac{1+\beta^{+}_i (g-1)}{1-\alpha^{+}_i (g-1)} \, 
\frac{1+\beta^{-}_i (g^{-1}-1)}
{1-\alpha^{-}_i (g^{-1}-1)}
\right]
\end{equation}
where
$$
0\le \alpha_i^\pm\,, \quad 0\le \beta_i^\pm \le 1\,, \quad 0\le \gamma^\pm 
$$
are parameters, 
are characters of $U(\infty)$ and conjectured that there
are no other characters. It was observed by Boyer \cite{Bo} and, independently, 
by Vershik and Kerov \cite{VK} that this conjecture is equivalent to the Schoenberg's
conjecture about the so-called
totally positive sequences (see below) which was already established by Edrei in \cite{Ed}
using some deep results about entire functions.

\subsection{} 

Vershik and Kerov also outlined a different and more direct proof which uses
approximation of characters of $U(\infty)$ by normalized
characters of $U(n)$. It follows from a general principle
due to Vershik (see \cite{V} and also \cite{Ol1}), 
that any character $\chi$ of $U(\infty)$ is a limit of a sequence
of normalized characters $\chi_n$ of $U(n)$ as $n\to\infty$ in the sense that
\begin{equation}\label{uniform}
\chi_n\Big|_{U(k)} \xrightarrow{\textup{  uniformly  }}  \chi\,\Big|_{U(k)}
\end{equation}
for any fixed $k=1,2,3,\dots$. In \cite{VK}, Vershik and Kerov gave necessary
and sufficient conditions for the convergence of $\{\chi_n\}$ and identified
the corresponding limits with functions \eqref{char}. 

This approximation
principle is a materialization of certain general ergodic theory 
ideas and is  closely akin to
some standard constructions in statistical physics such as
construction  of Gibbsian measures in an infinite volume by a thermodynamic
limit transition. In that case,
 one chooses a sequence of boxes which fill up the space (just as in
the above visualization of $U(\infty)$), for each box one picks some boundary condition
which specifies a Gibbsian measure (in our case, $\chi_n$), and one requires convergence of the
induced measures on all compact sets. 

\subsection{}

Note that the formula \eqref{char} is multiplicative in the eigenvalues of $g\in U(\infty)$.
This multiplicativity can be established a priori; as shown by Olshanski, see for example
\cite{Ol2}, such and more general 
multiplicativity are  very characteristic for representations of 
infinite-dimensional classical groups. 

Any character $\chi$ of $U(\infty)$ is therefore uniquely determined by its restriction
to $U(1)$
$$
x(z)=
\chi\left(\left[
\begin{matrix}
z \\
& 1 \\
&& \ddots
\end{matrix}
\right]
\right) = \sum_{k\in\Z} x_k\, z^k  \,.
$$
Conversely, any function $g\mapsto \det x(g)$ is a character of $U(\infty)$ provided
that it is positive definite, which means that its restriction to any $U(n)$ is a
nonnegative linear combination of the characters of $U(n)$, that is, of the rational 
Schur functions $s_\la$. 

The identity 
$$
\prod_{i=1}^n \left(\sum_{k\in\Z} x_k\, z_i^k \right)  = \sum_{\la=(\la_1\ge \dots \ge \la_n)\in\Z^n}  
\det\big[x_{\la_i-i+j}\big]_{i,j=1\dots n}\,
s_\la(z_1,\dots,z_n) \,,
$$
shows that this positivity is equivalent to
 the positivity of some (in fact, all) minors of the 
infinite Toeplitz matrix $\big[x_{j-i}\big]_{i,j\in\Z}$, which is precisely Schoenberg's
definition of a \emph{totally positive} sequence. 

In particular, the positivity of $2\times 2$ minors means that
$$
x_n^2 \ge x_{n-1} \, x_{n+1}  \,.
$$
Thus, one knows a priori that the restriction of any character of
$U(\infty)$ to $U(1)$ has log-concave multiplicities. 

\subsection{}\label{s24}

The question whether the same is true before the limit, that is, whether
the restriction of any irreducible representation of $U(n)$ to standard $U(1)$
has log-concave multiplicities, surfaced when we were working with G.~Olshanski
on a generalization of the Vershik--Kerov theorem \cite{VK}. Originally, this log-concavity
was needed to replace the uniform convergence in \eqref{uniform} by convergence
of Taylor series, see Section 3 in \cite{Ok1}. Eventually, 
in \cite{OO} it was replaced by a 
more elementary argument, but nonetheless this log-concavity is a valid 
question with interesting answer. 

As it turns out,  for any representation $V^\la$ of $U(n)$
the multiplicity of the irreducible representation $V^\mu$ of the standard
$U(k)\subset U(n)$ is a log-concave function of the pair 
$$
(\la,\mu)\in U(n)^\wedge \oplus U(k)^\wedge \,.
$$
In fact, one can say more. Without loss of generality, let us assume
that $\la$ is a partition and consider the space 
$$
V^{\la/\mu}=\Hom_{\,U(k)}(V^\mu\to V^\la)
$$ 
whose dimension is the multiplicity in question. The space $V^{\la/\mu}$
is an $U(n-k)$ module with character
given by the skew Schur function $s_{\la/\mu}$
$$
\tr_{\,V^{\la/\mu}}\left(\left[
\begin{matrix}
z_1\\
& z_2\\
&& \ddots
\end{matrix}
\right]
\right) = s_{\la/\mu} (z_1,z_2,\dots)\,.
$$
One has the following

\begin{theorem}[\cite{Ok1}]\label{t1} Suppose $(\la_i,\mu_i)$,
$i=1,2,3$, are partitions such that 
$$
(\la_2,\mu_2)=\tfrac12 (\la_1,\mu_1)+\tfrac12 (\la_3,\mu_3) \,.
$$
Then the following polynomial has nonnegative coefficients:  
\begin{equation}\label{slm}
s_{\la_2/\mu_2}^2 - s_{\la_1/\mu_1} s_{\la_3/\mu_3} \in \Z_{\ge 0}[z_1,z_2,\dots] \,.
\end{equation}
\end{theorem}

The coefficients of the polynomial $s_{\la/\mu}$ correspond to standard
tableaux of shape $\la/\mu$. In the proof of Theorem \ref{t1}, one 
constructs a certain transformation on pairs of standard tableaux and
proves that it is injective. 

Similar results for orthogonal and symplectic groups are also established in \cite{Ok1}.

\subsection{}

It is likely that \eqref{slm} is actually a nonnegative linear 
combination of Schur functions. One can propose a conjecture which
would, among other things, imply this property. 

Recall that the Littlewood-Richardson coefficients $c_{\la\mu\nu}$ are defined
by
$$
c_{\la\mu\nu}=\dim\left(V^\la\otimes V^\mu \otimes V^\nu\right)^{G}
$$ 
where the superscript $G$ stands for invariants of $G=U(n)$. If either of the
arguments of $c_{\la\mu\nu}$ is not a dominant weight, we set $c_{\la\mu\nu}=0$
by definition. Often, one uses the numbers 
$$
c^{\la}_{\mu\nu}=c_{\la^*\mu\nu}
$$
where $\la^*$ is the highest weight of the dual module $\left(V^\la\right)^*$
$$
(\la_1,\dots,\la_n)^*=(-\la_n,\dots,-\la_1) \,.
$$
The numbers $c^{\la}_{\mu\nu}$ are coefficients in the expansions
\begin{align*}
V^{\la/\mu}&=\sum_{\nu} c^{\la}_{\mu\nu} \, V^{\nu} \,, \\
V^{\mu}\otimes V^\nu &= \sum_{\la} c^{\la}_{\mu\nu} \, V^{\la} \,.
\end{align*}
\begin{conjecture}\label{c1} The function 
$$
(\la,\mu,\nu)\to \log\, c_{\la\mu\nu} 
$$
is concave.
\end{conjecture}

If true, this concavity would have some interesting applications. In particular, since
$$
c_{\la\mu\nu} = c_{\la\nu\mu}
$$
we conclude that
$$
c_{\la,\frac{\mu+\nu}2,\frac{\mu+\nu}2} \overset{?}\ge c_{\la\nu\mu} \,,
$$
provided $\frac{\mu+\nu}2$ is an integral weight. This is equivalent to the
inclusion of representations 
\begin{equation}\label{logV}
V^\nu\otimes V^\mu \overset{?}\subset \left(V^{\frac{\mu+\nu}2}\right)^{\otimes 2} \,, 
\end{equation}
which can be interpreted as saying that the representation valued function
\begin{equation}\label{Vla}
V:\la \mapsto V^\la
\end{equation}
is concave with respect to the natural ordering and tensor multiplication of
representations
\footnote{Remark that it follows from Weyl's dimension formula that
$$
\la \mapsto \log \dim V^\la
$$
is a concave function. That is, the function \eqref{Vla} considered as a function
into just vector spaces without group action is concave with respect to the tensor
product.}. 

If \eqref{logV} is true then we certainly have the following inclusion of $U(n)$-modules
\begin{equation}\label{vkl1}
\left(V^{\la_1}\otimes V^{\mu_1}\right) \otimes \left(V^{\la_3}\otimes V^{\mu_3}\right) 
\overset{?}\subset \left(V^{\la_2}\otimes V^{\mu_2}\right)^{\otimes 2} 
\end{equation}
for $(\la_i,\mu_i)$ as in Theorem \ref{t1}. 
The last inclusion is equivalent to
\begin{equation}\label{vkl2}
V^{\la_1/\mu_1}\otimes V^{\la_3/\mu_3}  \overset{?}\subset \left(V^{\la_2/\mu_2}\right)^{\otimes 2} \,. 
\end{equation}
Indeed, the equation \eqref{vkl2} is equivalent to
\begin{multline}\label{vkl3}
\left(\sum_{\nu_1} c_{\la_1^*\,\mu_1\,\nu_1} \, V^{\nu_1}\right) \otimes
\left(\sum_{\nu_3} c_{\la_3^*\,\mu_3\,\nu_3} \,V^{\nu_3}\right) \overset{?}\subset \\
\left(\sum_{\nu_2} c_{\la_2^*\,\mu_2\,\nu_2}\, V^{\nu_2}\right) \otimes
\left(\sum_{\nu_4} c_{\la_2^*\,\mu_2\,\nu_4}\, V^{\nu_4}\right)\,,
\end{multline}
whereas \eqref{vkl1} says that
\begin{multline}\label{vkl4}
\left(\sum_{\nu_1} c_{\la_1 \,\mu_1\,\nu_1^*}\, V^{\nu_1}\right) \otimes
\left(\sum_{\nu_3} c_{\la_3\, \mu_3\,\nu_3^*}\, V^{\nu_3}\right) \overset{?}\subset \\
\left(\sum_{\nu_2} c_{\la_2 \,\mu_2\,\nu_2^*}\, V^{\nu_2}\right) \otimes
\left(\sum_{\nu_4} c_{\la_2\, \mu_2\,\nu_4^*}\, V^{\nu_4}\right)\,. 
\end{multline}
To get \eqref{vkl3} from \eqref{vkl4}, take the dual space of everything, which
will replace $V^{\nu_i}$ by $V^{\nu_i^*}$, and then replace $\la_i$ by $\la_i^*$
and $\nu_i^*$ by $\nu_i$. The inclusion \eqref{vkl2}
is equivalent to Schur-positivity of \eqref{slm}. 

Similarly, the conjecture and the symmetry
$$
c_{\la\mu\nu}=c_{\nu\la\mu}
$$
imply  that 
$$
c_{\la',\mu',\nu'} \overset{?}\ge c_{\la\mu\nu}
$$
provided the weight
$$
\left(
\begin{matrix}
\la'\\ 
\mu'\\
\nu'
\end{matrix}
\right)
=
\left(
\begin{matrix}
\al & 0 & 1-\al \\
1-\al & \al & 0 \\
0 & 1-\al & \al
\end{matrix}
\right)
\left(
\begin{matrix}
\la\\ 
\mu\\
\nu
\end{matrix}
\right)\,, 
$$
is an integral weight and $0\le \al\le 1$. 

\subsection{}

Here is another implication of Conjecture \ref{c1} which is actually known to be true.
Concavity of $\log c_{\la\mu\nu}$ implies that the support 
$$
\supp c_{\la\mu\nu} =\{(\la,\mu,\nu),\,   c_{\la\mu\nu}\ne 0 \}
$$
is convex. In particular, since it contains the origin $(0,0,0)$, it is \emph{saturated},
meaning that
\begin{equation}\label{sat}
c_{k\la,k\mu,k\nu}\ne 0 \Rightarrow c_{\la,\mu,\nu}\ne 0 \,, 
\end{equation}
for any $k=2,3,\dots$. In fact, since $c_{0,0,0}=1$, Conjecture \ref{c1} implies that 
$$
c_{k\la,k\mu,k\nu} \overset{?}\le (c_{\la,\mu,\nu})^k \,.
$$
The saturation \eqref{sat} turns out to be a very important property, see \cite{F}. It has been recently
established by A.~Knutson and T.~Tao in \cite{KT}, see also \cite{Bu}.

\subsection{}
As already pointed out in Section \ref{s_irr}, log-concavity is (in contrast to so many
things in representation theory) not an additive property: it is totally destroyed 
by direct sums. 

It seems however likely that log-concavity should be a \emph{multiplicative} property,
that is, it should behave nicely with respect to tensor products. For example, recall that it is well known
and easy to prove that the convolution of two log-concave sequences is again log-concave.
This is equivalent to saying that the set of $U(1)$-modules with
log-concave multiplicities is closed under tensor products. 
 
This multiplicativity principle
fits together nicely with the above conjecture about tensor product multiplicities. 
 
\section{Log-concavity in the classical limit}\label{s3}

\subsection{}

Dealing with actual multiplicities may be a subtle business. 
Fortunately, many of these subtleties disappear in the classical
limit and much more general results can be established.

Let us assume that the phase space of our classical system is
an irreducible projective algebraic variety  $X\in\PP^N$ over $\C$ which is stable under
the action of a compact group $K\subset GL(N+1)$. We write $X$ in place of $M^{2n}$ 
to stress the fact that we are now working with projective
algebraic varieties which are allowed to be singular. 

Even for singular $X$, the moment map
$$
\phi: X \to \Lie(K)^*
$$
is still well defined as the restriction of the moment map for the 
$K$-action  on $\PP^N$. 
It is well known (see e.g.\ Theorem 6.5 in \cite{GS1}) that the function $\Omega(\xi)$ 
from Section \ref{s11}
describes the asymptotics of the multiplicities of $K$-modules in polynomials of 
very large degree on $X$. 

More concretely, let
$$
\C[X]=\oplus_{d=0}^\infty \C[X]_d
$$
be the homogeneous coordinate ring of $X$.  The space $\C[X]_d$ of degree $d$
polynomials on $X$ decomposes as a $K$-module
$$
\C[X]_d =\bigoplus_{\la\in K^\wedge} \Omega_d(\la)\, V^\la ,.
$$
One can view  $\Omega_k(\la)$ as a measure on $K^\wedge\subset\fh_+$. 
After proper normalization, the measures $\Omega_k(k\la)$ converge 
weakly to $\Omega(\xi)\, d\xi$ as $k\to\infty$  where $d\xi$ is the Lebesgue measure on $\fh_+$. 

In other words, for any $A\subset\fh_+$, the integral 
$\int_A\Omega(\xi)\,d\xi$ describes the leading asymptotics of the 
sum $\sum_{\la\in k A} \Omega_k(\la)$ as $k\to\infty$.  Hence, informally, $\Omega(\xi)$ is
 the multiplicity $\Omega_k(k\la)$ averaged over some infinitesimal neighborhood
of $\xi$. Such an averaging over infinitesimally close energy levels is a very natural thing to
do from the statistical physics perspective. 

\subsection{}
The function $\Omega(\xi)$ depends not only on the $K$-action on $X$ as such but also on
the embedding $X\subset \PP^N$,  where $\C^{N+1}$ is a representation space of $K$ or,
equivalently, of the complexification $G$ of $K$. 

In intrinsic terms, such an embedding is a  very ample invertible sheaf $L$ in the 
$G$-linearized Picard group $\Pic^G(X)$ of $X$. Write $\Omega(\xi,L)$ to
stress the dependence on both $\xi$ and $L$. Because $L$ enters the definition
of $\Omega(\xi,L)$ only via $L^{\otimes n}$, $n\to\infty$, the function 
$\Omega(\xi,L)$ is well-defined for any 
$$
L\in \Pic^G(X)\otimes_\Z\Q\,.
$$
Since ample sheafs form a semigroup
in $\Pic_G(X)$ it makes sense to ask whether $\log \Omega(\xi,L)$
is concave as a function of the pair $(\xi,L)$.  

In fact, we already saw
an example of such a bivariate concavity in Section \ref{s24} where the multiplicities
for restrictions from $U(n)$ to $U(k)$ turned out to be concave in the 
pair of highest weights.

\subsection{} 
In this setting, the log-concavity of $\Omega(\xi)$ and $\Omega(\xi,L)$ was
established in \cite{Ok2} and \cite{Ok3}, respectively, by using the classical
Brunn-Minkowski inequality of convex analysis. 

Here we want to use the same ideas to approach the problem from a slightly different angle.
Instead of looking at the weak limit of measures $\Omega_k(k\la)$, which involves averaging 
over infinitesimally close energy levels, we want to look at the asymptotics of
the sequence $\Omega_k(k\la)$ for some fixed $\la$. This is a natural thing to do
from the representation theory point of view. 

\subsection{} 

There is a standard trick which allows to dispose of the first variable $\xi$
in $\Omega(\xi,L)$ by enlarging the variety $X$. Indeed, by definition of the multiplicities
$\Omega_k$ we have
\begin{align*}
\Omega_k(k\la,L) &= \dim \left(H^0(X,L^{\otimes k}) \otimes V^{k\la^*}\right)^G\,, \\
&= \dim H^0 \left(X\times G/B,(L\boxtimes L_{\la^*})^{\otimes k}\right)^G\,, 
\end{align*}
where $\la^*\in K^\wedge$ is the highest weight of $(V^\la)^*$, $G/B=K/T$ is the 
flag variety of $K$, the sheaf $L_{\la^*}\in \Pic^G(G/B)$ corresponds to the map of $G/B$ onto
the orbit of the highest vector in $P(V^\la)$, and superscript $G$ denotes
$G$-invariants. So, without loss of generality,
we can assume that $\xi=0$.

\subsection{}

Recall that, by definition,  
$$
\C[X\q_LG] = \bigoplus_k H^0(X,L^{\otimes k})^G  
$$
is the homogeneous coordinate ring of the Geometric Invariant Theory quotient 
$X\q_LG$ corresponding to
$L\in\Pic^G(X)$. 
We write $X\q_LG$ in place of the standard $X\q G$ to stress the dependence on $L$.

The sequence 
\begin{equation}\label{dk}
\Omega_k(0,L)=\dim H^0(X,L^{\otimes k})^G\,, \quad k=0,1,2,\dots
\end{equation}
may fail to have a $k\to\infty$ asymptotics
for the following trivial reason. Consider the set
$$
\{k\,, H^0(X,L^{\otimes k})^G  \ne 0 \} \subset \Z_{\ge 0}\,.
$$
Since $X$ is irreducible, it 
is a semigroup and it either contains all sufficiently large integers
or lies in a proper subgroup if $\Z$. We want to avoid the latter case
because in that case the sequence \eqref{dk} does not have any asymptotics. So
we will replace $L$ by a suitable power of $L$ in that case.

By replacing $L$ by its power one can also achieve that $\C[X\q_LG]$ is generated
by its degree $1$ graded component and so we can assume this as well. Thus, we have
an embedding
$$
X\q_LG \subset \PP\left(\left(H^0(X,L)^G\right)^*\right) 
$$
and we denote by $\deg X\q_LG$ the degree of this embedding. It follows that
in this case
$$
\Omega_k(0,L) \sim  \deg X\q_LG \,\frac{k^{d}}{d!} \,, \quad  d=\dim X\q_LG\,,
$$
as $k$ goes to $\infty$. We now want to show that $\log\deg X\q_LG$ is a concave
function of $L$. 

\subsection{}
In fact, one can establish a more general fact. Let $Y$ be an irreducible
algebraic variety of dimension $d$ 
and let $\K=\C(Y)$ be the field of rational functions on
$Y$. Let $S\subset\K$ be a $\C$-linear subspace such that
$
1\in S
$
and which generates $\K$ as 
a field. The embedding $S\subset\K$ corresponds to a subvariety $Y_S\subset S^*$
which is birationally isomorphic to $Y$. Let $\deg Y_S$ denote its degree.

Given two such subspaces $S_1$ and $S_2$, denote by $S_1S_2$ the subspace generated
by all products $f_1 f_2$, where  $f_i\in S_i$. We will show that
\begin{equation}\label{BMdeg}
\sqrt[d]{\deg Y_{S_1 S_2}} \ge \sqrt[d]{\deg Y_{S_1}}  + \sqrt[d]{\deg Y_{S_2}}
\end{equation}
for any such pair $S_1$ and $S_2$. Since, clearly,
$$
\deg Y_{S^2} = 2^d\, \deg Y_S
$$
the inequality \eqref{BMdeg} implies that $\sqrt[d]{\deg Y_{S}}$, and
consequently, $\log\deg Y_{S}$ is a concave function of $S$. 

\subsection{}

In particular, \eqref{BMdeg} would imply the concavity of $\log\deg X\q_LG$. Indeed,
although the varieties $X\q_LG$ may not be isomorphic for different $L$,
they are always birationally isomorphic. Their common field of fractions
is the field $\K=\C(X)^G$ of rational $G$-invariants. 

Given some $L_1$ and $L_2$, pick some $\phi_i\in H^0(X,L_i)^G$. 
Replacing the $L_i$'s if necessary by their
multiples, we can assume
that $\C[X\q_{L_i}G]$ is generated by $S_i=\phi_i^{-1}H^0(X,L_i)^G\subset \K$. 
Since the algebra
$\C[X\q_{L_1\otimes L_2}G]$ contains the algebra generated by $S_1 S_2$, we get from
\eqref{BMdeg} the desired lower bound on the asymptotics of the dimensions of the
graded components of $\C[X\q_{L_1\otimes L_2}G]$.

\subsection{}

Now, in order to establish \eqref{BMdeg}, we will construct convex sets $\Delta_S\in\R^d$ 
of dimension $d=\dim Y$ such that 
$$
\deg Y_S  = d! \, \vol \Delta_S
$$
and
$$
\Delta_{S_1 S_2} \supset \Delta_{S_1} + \Delta_{S_2} \,.
$$
The inequality \eqref{BMdeg} will then follow immediately from the classical
Brunn-Minkowski inequality, see e.g.\ \cite{BZ}. See also e.g.\ the appendix
by A.~Khovanskii in \cite{BZ} for a discussion of the relationship between
classical inequalities of the convex analysis and algebraic geometry. For
example, the Alexandrov-Fenchel inequality, which is stronger than the
Brunn-Minkowski inequality, corresponds to the Hodge index theorem for surfaces. 

\subsection{}

The construction of $\Delta_S\in\R^d$ is similar to the definition of a Newton
polytope. 

Choose a smooth point $y\in Y$ which lies away from the 
singularities of the maps from $Y$ to $Y_{S_i}$ and their inverses. 
Choose a flag
of subvarieties
$$
Y \supset Y^1 \supset \dots \supset Y^d={y}\,, \quad \codim Y^k = k \,,
$$
which are all smooth at $y$. Fix some local equation $u_k$ of $Y^k$ in $Y^{k-1}$. 

This data give rise to a map
$$
\K\setminus 0 \owns f \mapsto \fv(f)=(\fv_1(f),\dots,\fv_d(f)) \in \Z^d \,, 
$$
where 
\begin{align*}
\fv_1(f)&=\ord_{Y^1} f\,,\\
\fv_2(f)&=\ord_{Y^2} 
\left. \left(f u_1^{-\fv_1(f)}\right) \right|_{Y^1}\,, \\
\fv_3(f)&= \ord_{Y^3} 
\left.\left(\left. \left(f u_1^{-\fv_1(f)}\right) \right|_{Y^1} u_2^{-\fv_2(f)}
\right)\right|_{Y^2}\dots
\end{align*}
and so on. It is clear that $\fv$ is a valuation, that is,
\begin{align}
\fv(fg)&=\fv(f)+\fv(g)\,,\label{val}\\
\fv(f+g)&\ge \min\{\fv(f),\fv(g)\}\,, \notag
\end{align}
where the ordering on $\Z^d$ is lexicographic. 

It is also clear that
the residue field of $\fv$ is isomorphic to $\C$ and hence for any
$\C$-linear subspace of $S\subset\K$ we have
\begin{equation}\label{dimv}
\dim_C S = |\fv(S\setminus 0)| \,.
\end{equation}

\subsection{}
By definition, set 
$$
\Gamma_S = \left\{(k,\fv(f)),f\in S^k\setminus 0\right\} \subset  \Z^{1+d} \,.
$$
It follows from \eqref{val} that $\Gamma_S$ is a semigroup. 

Denote by
$\Lambda_S\subset\Z^{1+d}$ the lattice generated by $\Gamma_S$.
Let $\nabla_S\subset\R^{1+d}$ be the closed convex cone generated by
$\Gamma_S$ and let $\Delta_S$ be the intersection of $\nabla_S$ with the 
subspace $(1,\R^d)\subset\R^{1+d}$. It is clear that
$$
\Delta_S = \overline{\left\{\frac{\fv(f)}k ,f\in S^k\setminus 0\right\}}\,,
$$
where bar denotes closure. 

Since the point $y$ corresponds to a smooth point of $Y_S$, there exist 
$$
f_0,f_1,\dots,f_d=1\in\C[Y_S]
$$
 such that
$$
\fv(f_k)=(\underbrace{0,\dots,0}_{\textup{$k$ times}}, 1,\dots) \,.
$$
Thus, $\Delta_S$ contains a $d$-dimensional simplex and so
\begin{equation}\label{dD}
\dim \Delta_S = d
\end{equation}

\subsection{} 
Let us now prove that 
$$
\vol \Delta_S = d! \, \deg Y_S\,,
$$
On the one hand, it follows that \eqref{dimv} and the definition of $\Gamma_S$ that
$$
\left|\Gamma_S\cap(k,\Z^d)\right| = \dim_C S^k 
\sim \deg Y_S \, \frac{k^d}{d!} \,, \quad k\to\infty \,.
$$
Since $\Gamma_S \subset \nabla_S\cap\Lambda_S$ we have
$$
\frac{\deg Y_S}{d!} \le \vol \Delta_S \,,
$$
where the normalization of Lebesgue measure is given by the intersection of 
the lattice $\Lambda_S$ with the 
subspace $(1,\R^d)\subset\R^{1+d}$. 

The inverse inequality will be deduced from the following
 results of Khovanskii \cite{Kh}.
Choose a sequence of finitely generated subgroups
$$
\Gamma_1 \subset \Gamma_2 \subset \dots\subset \Gamma_S
$$
such that $\Gamma_S = \bigcup
\Gamma_i$ and each $\Gamma_i$ generates the lattice $\Lambda_S$. Let
$\nabla_i$ denote the cone generated by $\Gamma_i$ and let $\Delta_i$
be the corresponding hyperplane section of $\nabla_i$. It is clear
that 
\begin{equation}\label{e33}
\Delta_1 \subset \Delta_2 \subset \dots  \subset\Delta_S=
\overline{\bigcup \Delta_i} \,.
\end{equation}
It is a theorem of Khovanskii, see Proposition 3 in Section 3 of \cite{Kh}, that
there exists vectors
$\gamma_i\in\Gamma_i$ such that
$$
(\nabla_i+\gamma_i)\cap\Lambda_S  \subset \Gamma_i \,.
$$
It follows that
$$
\frac{\left|\Gamma_i\cap(k,\Z^d)\right|}{k^d} \to \vol \Delta_i \,,  
\quad k\to\infty \,.
$$
Hence $\vol \Delta_i \le (d!)^{-1} \deg Y_S$ and \eqref{e33}
implies that
$$
\vol \Delta_S \le \frac{\deg Y_S}{d!} \,,
$$
as was to be shown. 

\subsection{}

The relation
\begin{equation}\label{inclD}
\Delta_{S_1 S_2} \supset \Delta_{S_1} + \Delta_{S_2} \,.
\end{equation}
follows immediately from \eqref{val} and the definition of $\Delta_S$. 
Now we are almost in position to finish the 
proof of \eqref{BMdeg} by applying the Brunn-Minkowski inequality.
One remaining detail is that our normalization of the volume $\vol \Delta_S$
depends on the lattice $\Lambda_S$. It is, however, clear that 
$$
\Lambda_{S_1 S_2} \subset \Lambda_{S_1}, \Lambda_{S_2} \,.
$$
Therefore, if we normalize the volume according to the lattice $\Lambda_{S_1 S_2}$ we
have
$$
\vol \Delta_{S_1 S_2} = \frac{\deg Y_{S_1 S_2}}{d!} \,, \quad 
\vol \Delta_{S_i} \ge  \frac{\deg Y_{S_i}}{d!} \,.
$$
This and the Brunn-Minkowski inequality applied to \eqref{inclD}
$$
\sqrt[d]{\vol \Delta_{S_1 S_2}} \ge \sqrt[d]{\vol \Delta_{S_1}}  + \sqrt[d]{\vol \Delta_{S_2}}
$$
completes the proof of \eqref{BMdeg}.

\end{document}